\documentclass{article}
\usepackage[utf8]{inputenc}

\usepackage{amsthm,amsfonts,amssymb,amsmath,epsf, verbatim}

\newtheorem{theorem}{Theorem}
\newtheorem{lemma}[theorem]{Lemma}
\newtheorem{corollary}[theorem]{Corollary}
\newtheorem{proposition}[theorem]{Proposition}

\newtheorem{conjecture}[theorem]{Conjecture}
\newtheorem{question}[theorem]{Question}

\title{Weighted Versions of the Arithmetic-Mean-Geometric Mean Inequality and Zaremba's Function}
\author{Tim McCormack\footnote{cortex@brainonfire.net}, Joshua Zelinsky\footnote{Hopkins School, New Haven, CT, USA,  zelinsky@gmail.com}}
\date{}

\begin{document}

\maketitle

\begin{abstract} We use the weighted version of the arithmetic-mean-geometric-mean inequality to motivate new results about Zaremba's function, $z(n) = \sum_{d|n} \frac{\log d}{d}$. We investigate record-setting values for $z(n)$ and the related function $v(n) = \frac{z(n)}{\log \tau(n)}$ where $\tau(n)$ is the number of divisors of $n$. We show that $v(n)$ takes on a maximum value and we give a list of all record-setting values for $v(n)$. Closely connected inequalities motivate the study of numbers which are pseudoperfect in a strong sense.
\end{abstract}

\section{Introduction}
Given a positive integer $n$, we will write $\sigma_k(n)$ as the sum of the $k$th powers of the divisors of $n$. We will write $\sigma(n)=\sigma_1(n)$ and write $\tau(n)$ for the number of divisors of $n$, noting that $\tau(n)=\sigma_0(n)$. Unless otherwise noted, when referring to `divisors' of a number, we will mean just the positive divisors.

Consideration of the means of the divisors of a number is an old problem. It is well known that if $n$ is a natural number, then the geometric mean of the divisors of $n$ is $\sqrt{n}$. The arithmetic mean of the set of divisors is just $\frac{\sigma(n)}{\tau(n)}$, and it is straightforward to see that the harmonic mean of the divisors of $n$ is 
$\frac{n\tau(n)}{\sigma(n)}$. These inequalities seem to have been first explicitly noted by Ore \cite{Ore1}. Ore also noted the following nice property: Let $D$ be the set of divisors of $n$, and let $h$ be the harmonic mean of $D$, and $a$ be the arithmetic mean of $D$. Then the  geometric mean of $h$ and $a$ is itself the geometric mean of $D$. This property does not hold in general if we replace $D$ with a general set of positive integers not arising from a set of divisors. Consideration of these means and when they take on integer values have been extensively studied. See for example, work by Bateman,  Erd\H{o}s, Pomerance and Straus \cite{BEPS}.  One might hope for something interesting from applying the arithmetic-mean-geometric-mean-harmonic-mean inequality (AM-GM-HM) to these divisor sets. Recall, the AM-GM-HM says that given a finite list of positive numbers, their arithmetic mean is at least their geometric mean which is at least their harmonic mean. Equality occurs exactly when all numbers in the list are equal. Applying this to the divisors of a number the bounds leads to the inequality: 

\begin{equation} \label{trivial application of AM-GM-HM}
 \frac{\sigma(n)}{\tau(n)} \geq \sqrt{n} \geq \frac{n\tau(n)}{\sigma(n)}.
\end{equation}

But Equation (\ref{trivial application of AM-GM-HM}) is essentially trivial and gives extremely weak bounds for all but very small $n$ since $\tau(n)$ grows slower than $n^{\epsilon}$ for any $\epsilon >0$ . At some level, the triviality should not be surprising; these inequalities hold for any positive set of numbers, and the inequalities do not ``know'' that we are dealing with sets of natural numbers, much less sets of divisors. So the room to get number theoretic content out of these is small.

However, there is a weighted version of the AM-GM-HM. A question then arises whether by using carefully chosen weights we can get statements with some number theoretic content. The main goal of this paper is to show that the answer is ``yes, but only weakly.'' \footnote{While this paper was under review, a paper by S\'{a}ndor appeared \cite{Sandor} which also looks at using the weighted weighted version of the AM-GM to derive some similar inequalities.}

Recall, the weighted arithmetic-geometric mean inequality says that 
\begin{lemma} If $\alpha_1, \alpha_2 \cdots ,\alpha_k$ are real numbers, with $\sum_{i=1}^k \alpha_i=1$, and $x_1, x_2, \ldots ,x_k$ are positive real numbers, then \label{Weighted AM GM}
$$\prod x_i^{\alpha_i} \leq \sum \alpha_i x_i.$$
\end{lemma}

In what follows, we will refer to Lemma \ref{Weighted AM GM} as simply the weighted AM-GM inequality. 

The goal of this paper is to use the weighted AM-GM inequality, and various tightened versions of it, to obtain inequalities which while, not strong, are to our knowledge not in the literature. Some of our bounds with the weighted AM-GM inequality will be surpassed by other bounds in papers which use other methods (including parts of this one), and thus one way of thinking about the use of the weighted inequality here is as a motivator for examining the relationship between different functions. In the course of doing so, we will also introduce a new class of pseudoperfect numbers which will match almost precisely the broadest context where some of our proofs apply.

Many of our results will concern  Zaremba's function, which is defined as $z(n) = \sum_{d|n}\frac{\log d}{d}$. Notice that $z(n)= \log \prod_{d|n}d^{\frac{1}{d}}$.  This function was first introduced by Stanis\l{}aw Zaremba who used it to study certain lattices relevant for numerical integration \cite{Zaremba}. Subsequently, the function was investigated by Erd\H{o}s and Zaremba \cite{EZ}, who proved that

$$\limsup \frac{z(n)}{(\log \log n)^2} = e^\gamma,$$ where $\gamma$ is Euler's constant.

The average value of $z(n)$ does not appear to be explicitly noted in the literature, but it is straightforward to find. 
\begin{proposition} We  have
$$\sum_{n \leq x}z(n) = x(-\zeta'(2)) + O(\log^2 x).$$
\end{proposition}
\begin{proof} We have \begin{equation}\sum_{n \leq x} z(n) = \sum_{d \leq x}\left\lfloor \frac{x}{d}\right\rfloor \frac{\log d}{d} = \sum_{d \leq x}\left( \frac{x}{d} +O(1)\right) \frac{\log d}{d} = \sum_{d \leq x}\frac{\log d}{d^2} +O\left(\sum_{d \leq x} \frac{\log d}{d} \right). \end{equation}

We now note that $$\sum_{d \leq x} \frac{\log d}{d}= O ((\log x)^2), $$

and that $$x\sum_{d \leq x}\frac{\log d}{d^2} = x\left(\sum_{d=1}^{\infty} \frac{\log d}{d^2}+O\left(\frac{\log x}{x}\right)\right) = x(-\zeta'(2)) + O(\log x).\qedhere$$ 
\end{proof}

Weber \cite{Weber} is also notable for examining a variety of functions whose behavior is similar but not exactly the same as $z(n)$, such as replacing the $\log d$ in the definition with $(\log d)(\log \log d)$. 

We also will use various upper bounds on $z(n)$ as functions of which primes divide $n$. 

\begin{lemma}\label{Lemma 4.3 from Weber}(Lemma 4.3 from \cite{Weber}) For all $n$, $$z(n) \leq\left( \prod_{p|n} \frac{p}{p-1} \right)\left(\sum_{p|n} \frac{\log p}{p-1}\right). $$ Here $p$ ranges only over primes.
\end{lemma}

Lemma \ref{Lemma 4.3 from Weber} is in fact a corollary of another Lemma from that paper, which will be used frequently in this paper. 

\begin{lemma} If $n=p_1^{a_1}p_2^{a_2} \cdots p_k^{a_k}$, then 
\begin{equation}\label{Erdos-Zaremba rearrangement} z(n) = \sum_{i=1}^k \sum_{j=1}^{\alpha_i} \frac{j \log p_i}{p_i^j}h\left(\frac{n}{p_i^{\alpha_i}}\right). \end{equation} 
\end{lemma}

Note that a similar version of this Lemma in Erd\H{o}s and Zaremba is incorrect. This is the corrected version in Weber. They have a $h(\frac{n}{p^j})$ term when it should be $h(\frac{n}{p^{\alpha_i}})$.

We also have 
\begin{lemma} Assume that $n=p_1^{a_1}p_2^{a_2} \cdots p_s^{a_s}q_1q_2\cdots q_t$ where the $p_i$ and $q_i$ are all distinct primes. Then $$z(n) \leq \left(\sum_{i=1}^s \frac{\log p_i}{p_i -1} + \sum_{i=1}^t \frac{\log q_i}{q_i}\right)\left(\prod_{i=1}^s \frac{p_i}{p_i-1}\right)\left(\prod_{i=1}^s \frac{q_i+1}{q_i}\right). $$ \label{Tightened Weber 4.3} 
\end{lemma}
\begin{proof} A version of this also proven in Weber, we include a proof here. Assume we have such an $n$. We set $P=p_1^{a_1}p_2^{a_2} \cdots p_s^{a_s}$, and $Q=q_1q_2\cdots q_t$. Note that $(P,Q)=1$, and that $n=PQ$. Then, from Lemma \ref{Erdos-Zaremba rearrangement}, 
we have \begin{equation} \label{Multistage inequality}
\begin{split}
z(n) & =  \sum_{i=1}^s \sum_{j=1}^{\alpha_i} \frac{j \log p_i}{p_i^j}h\left(\frac{PQ}{p_i^{\alpha_i}}\right) + \sum_{i=1}^t \frac{\log q_i}{q_i}h\left(\frac{PQ}{q_i}\right) \\
 & \leq h(PQ)\left(\sum_{i=1}^s \sum_{j=1}^{\alpha_i} \frac{j \log p_i}{p_i^j} + \sum_{i=1}^t \frac{\log q_i}{q_i}\right). 
\end{split}
\end{equation}

We have that $$h(PQ) \leq \prod_{i=1}^s \frac{p_i}{p_i-1} \prod_{i=1}^t \frac{q_i}{q_i+1},$$

which combines with $\sum_{j=1}^{\alpha_i} \frac{j \log p_i}{p_i^j} \leq \frac{ \log p_i}{p_i-1}$, to give our desired inequality.
\end{proof}

Note that when $Q=1$ (equivalently $t=0$), Lemma \ref{Tightened Weber 4.3} becomes just Lemma \ref{Lemma 4.3 from Weber}, that is Weber's 4.3 Lemma. 

In what follows, we will also write $h(n)=\sigma_{-1}(n) = \sum_{d|n}\frac{1}{d}$, and note that $h(n)= \frac{\sigma(n)}{n}$.  Recall that a number is said to be perfect if $\sigma(n)=2n$, or equivalently that $h(n)=2$. A number is said to be abundant if $\sigma(n) > 2n$, and deficient if $\sigma(n) <2n$.

We will write $H(n)= \prod_{p|n} \frac{p}{p-1}$. Note that $H(n)= \frac{n}{\phi(n)}$ where $\phi(n)$ is the Euler totient function.  We will set $M(x) = \prod_{p\leq x} \frac{p}{p-1}$. Also, given a positive integer $n$, with $t$ positive divisors, we will write the divisors of $n$ in increasing order as $d_0$, $d_1, \cdots d_{t-1}$. Note that we have started our index at $0$ and ended at $t-1$. This notation will be convenient since often we will want to take a sum which runs of over all divisors of an integer which are greater than 1, rather than a sum over all positive divisors.

We recall here a few well-known facts which will be useful to us. In particular

\begin{lemma} $h(n) \leq H(n)$ with equality if and only if $n=1$.\label{h(n) < H(n)} And for any $n$, $$\lim_{k \rightarrow \infty}h(n^k) = H(n).$$
\end{lemma}

\begin{lemma}(Mertens' theorems)\label{Mertens'theorem} We have  $$M(x) \sim e^\gamma \log x$$ where $\gamma$ is Euler's constant. We also have

$$\sum_{p \leq x} \frac{\log p}{p} \sim \log x,$$ where the sum on the right-hand side is over all primes at  most $x$.
\end{lemma}

Lemma \ref{Mertens'theorem} can be found in  most standard number theory textbooks. Stronger versions of both equations in Lemma \ref{Mertens'theorem} with explicit error term can be found in a variety of papers such as Rosser and Schoenfeld's classic paper \cite{RosserSchoenfeld1962}.

\begin{lemma}(Euclid-Euler Theorem). \label{Euclid-Euler)} If $n$ is even, then $n$ is perfect if and only if $n= (2^p-1)(2^{p-1})$ where $2^p-1$ is prime.
\end{lemma}

Note that if $2^p-1$ is prime, then $p$ is prime, but the converse is not true.  This fact motivates our choice of the letter $p$ in the above lemma. 








\section{Upper bounds for $z(n)$ in terms of $\tau(n)$}
One gets a particularly nice result from applying the weighted version of the AM-GM inequality to perfect numbers. Let $n$ be a perfect number. We will write $d_1, d_2 \cdots d_k$ as the positive divisors of $n$ greater than 1 in increasing order. We may apply the weighted AM-GM inequality with $\alpha_i = \frac{1}{d_i}$, and $x_i =d_i$.  Note that $\sum_{i=1}^k \alpha_i=1$ and so we may apply the weighted AM-GM inequality using the $a_i$ as our weights. (This is why we need to often throw 1 out as one of our divisors, because otherwise our weights will sum to 2.)

We then apply the weighted AM-GM inequality to obtain that

\begin{theorem} If $n$ is a perfect number, then
$$\prod_{d|n} d^{\frac{1}{d}} < \tau(n)-1.$$
Equivalently, 
 $z(n) < \log (\tau(n)-1).$ \label{tau lower bound for perfect form}
\end{theorem}

Note that this inequality is false for some non-perfect $n$. For example, it fails whenever $n$ is prime, and it also fails for $n=24$. Thus, $n$ being perfect is really necessary in Theorem \ref{tau lower bound for perfect form}.  Note also that when $n$ is perfect, $z(n)$ is the Shannon entropy of the probability distribution where each divisor $d>1$ has probability $\frac{1}{d}$.

There's still not so much number theoretic content in this form, since we're just using that the sum of $\frac{1}{d_i}$ has the correct sum. In fact, we can make a version of this for deficient numbers which is almost as tight. Let $h(n)=\frac{\sigma(n)}{n}.$ Then as before we set for each $i \geq 1$, $\alpha_i=\frac{1}{d_i}$. We set $x_i = d_i$, for all $0 \leq i \leq t $ (including 1). However, we set $\alpha_0 =2 - h(n)$. Notice that we do in fact have $\sum \alpha_i =1$, with $\alpha_0$ essentially making up for the deficiency. We then apply the weighted AM-GM inequality, note that $1$ to any power is 1 and obtain the following.

\begin{theorem} If $n$ is a deficient number then 

$$\prod_{d|n} d^{\frac{1}{d}} \leq \tau(n) + 1 - h(n).$$

Equivalently,

$$z(n) \leq \log \left(\tau(n) + 1 - h(n) \right) .$$
\end{theorem}

But even this might overstate the number-theoretic content. We could have chosen the $x_i$ to be any set of real numbers all greater than 1 whose sum of reciprocals is at most 1 and gotten an essentially similar result. 

In fact, the situation is slightly worse than that. Recall, a number $n$ is said to be pseudoperfect if the sum of some subset of its proper divisors  adds up to $n$.  Trivially, any perfect number is also pseudoperfect. Theorem \ref{tau lower bound for perfect form} applies just as well to pseudoperfect numbers where we replace $\tau(n)-1$ with the number of elements in the relevant divisor sum. In what follows, most but not all of our results about perfect numbers will apply to pseudoperfect numbers. 

Given a positive integer $n$, when $n$ is pseudoperfect, we will write $S_0$  to be a subset of its positive proper divisors which sum to $n$. We will also write $S_1(n)$ to be a subset of all the divisors of $n$ (possibly including $n$), and with the sum of $S_1$ adding to $2n$. Notice that any $S_0$ gives rise to a valid $S_1$ by $S_1 \cup \{n\}$.)  We will refer to $S_0$  as a divisor set of $n$. Notice that we refer to {\emph{a}} divisor set of $n$, rather than {\emph{the}} divisor set of $n$, since there may be more than one. In particular, if $n_1$ and $n_2$ are distinct pseudoperfect numbers, then $n_1n_2$ is pseudoperfect with at least two different choices of $S_0$ one inherited from $n_1$ and one inherited from $n_2$. For example, 6 is pseudoperfect, since $1+2+3=6$, and $20$ is pseudoperfect with $1+4+5+10=20$, We have  $$1(20)+2(20)+3(20)=6(20)=6(1)+6(4)+6(5)+6(10).$$
Thus, we have as choices for $S_0$ both $\{20,40,60\}$ as well as $\{6, 24, 30, 60\}$. (Some pseudoperfect numbers have more than one choice of $S_0$ even though they are not inheriting them from another set.)

With our new notation to discuss not just perfect but also pseudoperfect numbers,  Theorem \ref{tau lower bound for perfect form}  becomes the more general claim:

\begin{theorem} Let $n$ be a pseudoperfect number with divisor set $S_0$. Then $$\sum_{d \in S_0} \frac{\log d}{d} \leq  \log |S_0|. $$
\end{theorem}

This result is also not only true for all pseudoperfect numbers, but it also is very weak even for even perfect numbers which we know actually exist. In particular, a consequence of the Euclid-Euler theorem on even perfect numbers is that if $n$ is an even perfect number, then $$z(n) < \left(\sum_{i=1}^\infty \frac{\log (2^i)}{2^i}\right)\left(\sum_{i=1}^\infty \frac{\log (3^i)}{3^i}\right) = \frac{3}{2}(\log 2)(\log 3).$$
In fact, if there are infinitely many even perfect numbers and $P_i$ is the $i$th even perfect number then $$\lim_{i \rightarrow \infty} z(P_i) = 2\log 2.$$

So, at least in their case, this is a very weak bound. It seems like getting stronger lower and upper bounds on $z(n)$ when $n$ is an odd perfect number would be worthwhile.

One may wonder also how tightly one can bound $z(n)$ from above in terms of $\tau(n)$ when there are no special restrictions on $n$. We first have the following obvious bound.

\begin{proposition} We have $z(n) = O\left(\log (\tau(n))^2\right )$. \label{z(n) upper bound or arbitrary integers}
\end{proposition}
\begin{proof} Note that $\frac{\log x}{x}$ is a decreasing function for $x \geq e$, and thus we have $z(n) \leq \sum_{i=1}^{\tau(n)+1} \frac{\log i}{i} = O((\log \tau(n))^2).$
\end{proof}

However, we can with a little work get a better result, and again just using the weighted AM-GM inequality.

\begin{lemma} Let $n$ be a positive integer greater than 1. Then $$z(n) < \left(h(n)-1\right)\log \frac{\tau(n)-1}{h(n)-1}.$$ \label{z(n) < h log tau}
\end{lemma}

\begin{proof} This is essentially just applying the weighted AM-GM inequality with $x_i$ running over $d|n$ for $d>1$ and 
taking as weights $\alpha_i = \frac{\frac{1}{x_i}}{\sum_{i=1}^{\tau(n)-1} \frac{1}{x_i}}.$
\end{proof}
Note also that the inequality in Lemma \ref{z(n) < h log tau} is strongest when $h(n)-1$ is very small or very large. This becomes essentially our earlier inequality for perfect numbers when $h(n)=2$. We may use Lemma \ref{z(n) < h log tau}  to get a tighter bound for $z(n)$ purely in terms of $\tau(n)$. 

\begin{proposition}
 $z(n)= O\left( (\log \tau(n)(\log \log \tau(n) \right)$.
\end{proposition}
\begin{proof} By Lemma \ref{z(n) < h log tau}, we need only show that $h(n) = O( \log \log \tau(n))$. To see this, note that $n$ has at most $\log_2 \tau(n)$ distinct prime divisors. Thus, if $P_j$ is the $j$th prime number then we have $$h(n) < H(n) \leq \prod_{p \leq P_{\lfloor \log_2 \tau(n) \rfloor}} \frac{p}{p-1} .$$ We have  $$ \prod_{p \leq P_{\lfloor \log_2 \tau(n) \rfloor}} \frac{p}{p-1}  = M(P_{\lfloor \log_2 \tau(n) \rfloor}) = O\left(\log (\log_2 \tau(n) \log \log \tau(n))\right) = O(\log \log \tau(n)).  $$
where the penultimate step is an application of the Prime Number Theorem together with Lemma \ref{Mertens'theorem}.
\end{proof}

Comparing the size of $z(n)$ to $\tau(n)$ motivates the definitions of the following functions. 
Define $v(n)$ by $$v(n) = \frac{z(n)}{\log \tau(n)}.$$

Two natural sequences of integers are the record setters for $v(n)$ and $z(n)$. By the record setters of a function $f(n)$ we mean the set $R_f$ defined by 
$$R_f:=\{n \in \mathbb{N}| \forall m \in \mathbb{N}, m < n \implies f(m) < f(n)\}.$$

That is, a number is a record setter for a function $f(x)$ if $f(n)$ is larger than $f(i)$ for any $0 < i<n$.
Record setters for arithmetic functions have been studied since Hardy and Ramanujan. Examples include $R_\tau$, called the highly composite numbers, $R_{h}$, called the superabundant numbers, and $R_\sigma$, called the highly abundant numbers.\cite{Alaoglu and Erdos}\cite{Ramanujan HC}

Neither $R_z$ or $R_v$ is currently in the OEIS.
For $v(n)$, the first few record setters are
2, 3, 4, 6, 12, 24, 36, 48, 60, 120, 180, 240, 360, 720, 840, 1260, 1680, 2520, 5040, 10080, 15120.

One may wonder if $v(n)$ and $z(n)$ ever have different record setters aside from the initial 1. And in fact, they do start off the exact same way. The first few record setters for $z(n)$ are also 
1, 2, 3, 4, 6, 12, 24, 36, 48, 60, 120, 180, 240, 360, 720, 840,	1260,	 1680, 2520, 5040, 10080, 15120,

However, they disagree at the next term, 25200, which is a record setter for $z(n)$ but not for $v(n).$ More terms for both are included in the Appendix. If one looks at that data, one will notice a few patterns. In that list, whenever $v(n)$ has a record setter, so does $z(n)$. We will show this is the case by showing that in fact the set of record setters for $v(n)$ is a finite list where we explicitly give all of them.  It  turns out that the very last record setter for $v(n)$ is $321253732800$. 

Many consecutive elements in  $R_z$ have the same number of divisors. Can that happen infinitely often? One might think that whenever this happens, that $v(n)$ is now also a record setter. But this can fail if $\tau(n)$ is already too large in that range to contribute to a new record setter. For example, 4497552259200 and 4658179125600 are both record setters for $z(n)$ and both have 9216 divisors. However, the second is not a record setter for $v(n)$, since although $v(4658179125600) > v(4497552259200)$, both have smaller $v(n)$ values than the previous record setter.

The highly composite numbers, superabundant numbers, and highly abundant numbers share the property that all of them are ``waterfall numbers.'' (Sometimes called ``Hardy-Ramanujan'' numbers.)  These numbers are $n$ of the form $n=p_1^{a_1}p_2^{a_2}\cdots p_k^{a_k}$ where $p_i$ is the $i$th prime, and where $a_{b} \geq a_c$ whenever $b < c$. These numbers form OEIS sequence A025487.\cite{OEIS A025487}

We have a similar property for $R_z$ and $R_v$. 

\begin{proposition} With the exception of 3, $R_z$ and $R_v$ consist only of waterfall numbers.
\end{proposition}
\begin{proof}
We sketch the proof for $R_z$. The proof for $v(n)$ is nearly identical. The same basic proof as for the other record setting sets works; if a number $n$ is a record-setter and one has $a_{b} < a_c$ for some $b < c$, then one could swap $a_b$ and $a_c$ to get a smaller number $n'$ where $z(n') > z(n)$, and hence $n$ would not be a record setter. The only difficulty for $n=3$ arises because $\frac{\log x}{x}$ is not a decreasing function until $x \geq e$, and in fact $\frac{\log 2}{2} < \frac{\log 3}{3}$. To handle this case, note that one needs only consider possible record setters $n$  where $3||n$ or $n=3^k$.  If $n=3^k$ then it can only be a record setter when $k=0$ or $k=1$. Thus, we may assume that $n$ is either divisible by 2 or divisible by $p$ for some $p>3$. If $n$ is not divisible by $2$, then we may  take some prime $p >3 $ replace all copies of that prime $p$ in the prime factorization of $n$ with $2$s, and get a smaller number with a larger $z$ value. Thus, we may assume that $n$ is divisible by some power of 2 and a single 3, and then some larger primes, possibly raised to larger powers. But then the argument that $\frac{\log x}{x}$  is a decreasing function may then again be applied. 
\end{proof}

Since every element of $R_v$ and $R_z$ is either 3 or a waterfall number, finding all record setters less than a given $x$ is straightforward. To do this efficiently, one use the fact that $n$ being a waterfall numbers has the equivalent formulation that $n$ can be represented as the products of (not necessarily distinct) primorials. 


It is not hard to show that $h(n)$ grows arbitrarily large as $n$ ranges over elements of $R_z$. If $h(n)$ is bounded, then the behavior is different. In particular,

\begin{conjecture} Fix a constant $C>1$, and let $n_1, n_2 \cdots$ be an increasing sequence of positive integers satisfying the following conditions.\label{z small compared to tau when h fixed}
\begin{enumerate}
    \item For all $i$, $h(n_i) \leq C$.
    \item $\lim_{i \rightarrow \infty} \tau(n_i)=\infty$.
\end{enumerate}
 Then $$z(n_i) = o(\log \tau(n)).$$
\end{conjecture}

We will prove Conjecture \ref{z small compared to tau when h fixed} in some limited cases.

\begin{proposition} Assume that there are infinitely many even perfect numbers. Then the sequence of even perfect numbers satisfies Conjecture \ref{z small compared to tau when h fixed}. \label{Even perfects satisfy z small compared to tau}
\end{proposition}
\begin{proof} Let $n_i$ be the $i$th even perfect number, and write $n_i=(2^{p_i-1})(2^{p_i}-1)$ which we can do by Lemma \ref{Euclid-Euler)}. We have  that $\lim_{i \rightarrow \infty} \tau(n_i)=\infty$, so we need only prove that $z(n_i)$ is bounded above by some constant. By Lemma \ref{Lemma 4.3 from Weber},
$$z(n_i) \leq \frac{2}{2-1}\frac{(2^{p_i}-1)}{(2^{p_i}-2)}\left(\frac{\log 2}{2-1} + \frac{\log(2^{p_i}-1))}{2^{p_i}-2}\right).$$
Since $\frac{\log x}{x-1}$, and $\frac{x}{x-1}$ are both decreasing functions, we have that 
$$\frac{2}{2-1}\frac{(2^{p_i}-1)}{(2^{p_i}-2)}\left(\frac{\log 2}{2-1} + \frac{\log(2^{p_i}-1))}{2^{p_i}-2}\right) \leq 2\left(\frac{3}{2}\right)\left(\log 2 + \frac{\log 3}{3}\right),$$ and so we are done.
\end{proof}

Proposition \ref{Even perfects satisfy z small compared to tau} required the assumption that there are infinitely many even perfect numbers. We can make a version of the statement which applies to a broader class of numbers which does not rely on any unproven assumption. In particular, a number is said to be primitive non-deficient if it satisfies $h(n) \geq 2n$, and $h(m)< h(2m)$ for any $m$ where $m|n$, and $m<n$. For any odd prime $p$, any number $n$ of the form $$n=p2^{\lfloor \log_2 p \rfloor}$$ is primitive non-deficient. Essentially, we take $p$ and we multiply it by the smallest power of 2 which is at least $\frac{p}{2}$. Even perfect numbers are of this form, and are essentially the case where the power of 2 is just barely above $\frac{p}{2}$. Numbers of this form also all satisfy $h(n) \leq 3$.

We may then modify the argument from Proposition \ref{Even perfects satisfy z small compared to tau} to handle this specific class of primitive non-deficient numbers.

\begin{proposition} Let $p_i$ be the $i$th odd prime. Then the sequence $s_i =p_i2^{\lfloor \log_2 p_i \rfloor} $ satisfies Conjecture \ref{z small compared to tau when h fixed}.
\end{proposition}

In fact, for the sequence $s_i$, we have the stronger result that $\lim_{i \rightarrow \infty} s_i = 4 \log 2$.

We have one other sequence of primitive non-deficient numbers which satisfy Conjecture \ref{z small compared to tau when h fixed}. 

\begin{proposition} Let $p_i$ be the $i$th prime number and let $ c_j=p_jp_{j+1}\cdots p_{k_j}$ where $k_j$ is the smallest $k$ such that the product is not deficient. Then $c_j$ satisfies Conjecture \ref{z small compared to tau when h fixed}. In particular, we have that $\log \tau(c_j)$ is asymptotic to $(\log 2)\frac{p_j^2}{2 \log p_j}$. We also have that
$j$, $z(c_j)=O(\log p_j)$.\label{aj has z not much more than log pj}
\end{proposition}
\begin{proof} It follows from the techniques used in Section 5 of \cite{Zelinskybig}, that for large $j$, $k_j$ is asymptotic to $\frac{p_j^2}{2 \log p_j}$, so $\log \tau(c_j)$ is asymptotic to $(\log 2)\frac{p_j^2}{2 \log p_j}$.

That we may bound $z(c_j)$ above then follows from applying Lemma \ref{Lemma 4.3 from Weber}, and using that $\sum_{p \leq x}\frac{\log p}{p} = O(\log x)$.  In particular, $z_(c_j) = O(\log p_{k_j})$, but since $p_{k_j} \sim \frac{p_j^2}{2}$ the bound follows.  
\end{proof}

We now turn our attention to verifying that $R_v$ is finite, and that $R_v$'s largest element is $321253732800$. 

\begin{lemma} If $n$ is a record setter for $z(n)$ or $v(n)$, and $13|n$, then $4|n$.
\end{lemma}
\begin{proof} 
Set $m=3(5)\cdots p_{k-1}$ where $p_{k-1} \geq 7$. Set
$n_1= 4m$, and $n_2=2mp_k$. We have that $$\tau(n_1) = 3 (2^{k-2}) < 2^k = \tau(n_2),$$ so the Lemma will be proved if we can show that
$z(n_1) > z(n_2).$ This is sufficient because record setters must be waterfall numbers, so if $a_1=1$, all the other prime factors must be raised to the first also.

Our plan is to use Equation (\ref{Erdos-Zaremba rearrangement})  to directly estimate $z(n_1)-z(n_2)$ by splitting it up into two terms. One term has the contribution from both $i=1$ $i=k$, and the other term has the contribution for $i$ ranging from $2$ to $k-1$. We will show that both are positive. We have 

$$z(n_1) - z(n_2) = A+B$$
where $$A= \frac{\log 2}{2}h\left(\frac{n_1}{4}\right) + \frac{2\log2}{4}h\left(\frac{n_1}{4}\right) - \frac{\log 2}{2}h\left(\frac{n_2}{2}\right) - \frac{\log p_k}{p_k}h\left(\frac{n_2}{p_k}\right),  $$

and $$B= \sum_{i=2}^{k-1} \frac{\log p_i}{p_i}\left(h\left(\frac{n_1}{p_i}\right) -h\left(\frac{n_2}{p_i}\right)\right). $$

Since $h(n)$ is a multiplicative function we may rewrite $A$ as 
$$A=h(m)\left(\log 2 - \frac{\log 2}{2} \frac{p_k+1}{p_k} - \frac{3}{2} \frac{\log p_k}{p_k} \right),$$

and the quantity in the parentheses is positive as long as $p_k \geq 13$.

For $B$, we will show that for any $i$, where $2 \leq i \leq k-1$,
that $$h\left(\frac{n_1}{p_i}\right) - h\left(\frac{n_2}{p_i}\right) \geq 0.$$

Since $n_1 = 4m$, and $n_2 = 2mp_k$, this is the same as verifying that 
$$\frac{7}{4}h\left(\frac{m}{p_i}\right) - \frac{3}{2}\left(\frac{p_k+1}{p_k}\right)h\left(\frac{m}{p_i}\right) \geq 0.$$ This last inequality is easily seen to be true when $p_k \geq 7$.

\end{proof}

We can generalize this lemma to get the following:

\begin{lemma} Let $p_i$ be the $i$th prime number, and fix some $\ell$. Assume that $n$ is a record setter for $v(n)$ or $z(n)$  and $p_k$ satisfies \label{repeated prime factor lemma}

\begin{enumerate} 
\item $p_k > p_{\ell}(p_\ell +1)$
\item $\log p_\ell + \frac{(2\log p_\ell) }{p_\ell} \geq  (\log p_\ell)\frac{p_k+1}{p_k} + \frac{\log p_k}{p_k}(p_\ell +1)  $ 
\item $p_k|n$.
\end{enumerate}
then  $p_\ell^2|n$. 

\end{lemma}
\begin{proof} Assume that $p_k$ satisfies the above three conditions. Note that the first condition implies $p_k > p_\ell$. Set \begin{equation}m= 2^{a_1} \cdots p_{\ell-1}^{a_{\ell-1}}p_{\ell+1}\cdots p_{\ell-1}. \end{equation}
Notice that $m$ does not include a $p_\ell$ term, and  . We will set 
$n_1 = p_\ell^2m$, and set $n_2 = p_\ell mp_k$. As with the previous Lemma, the result will be proven if we can show that $z(n_1) \geq z(n_2)$. (As before we get that $\tau(n_1) < \tau(n_2)$ by the same logic.)

We apply Lemma \ref{Erdos-Zaremba rearrangement} to get that 
$$z(n_1)-z(n_2) = A+B,$$ where $A$ has the terms arising from $i=\ell$ and $i=k$, and $B$ has all the other terms.

Then $$A=\frac{\log p_\ell}{p_\ell}h\left(\frac{n_1}{p_\ell^2} \right) + \frac{2\log p_\ell}{p_\ell^2}h\left(\frac{n_1}{p_\ell^2}\right) - \frac{\log p_\ell}{p_\ell} h\left(\frac{n_2}{p_\ell}\right) -\frac{\log p_k}{p_k}h\left(\frac{n_2}{p_k}\right) .$$

Applying that $h(n)$ is multiplicative, we may rewrite this as 
$$A=\frac{h(m)}{p_\ell}\left(\log p_\ell + \frac{2\log p_\ell }{p_\ell} - (\log p_\ell)\frac{p_k+1}{p_k} - \frac{\log p_k}{p_k}(p_\ell +1) \right). $$

The quantity inside the parentheses is positive by condition 2 above.

We have $$B=\sum_{i=1, i \neq \ell }^{k-1} \frac{\log p_i}{p_i}\left(h\left(\frac{n_1}{p_i}\right) -h\left(\frac{n_2}{p_i}\right)\right). $$

So we just need to verify that for $i$ with $2 \leq i \leq k-1$, we have $$h\left(\frac{n_1}{p_i}\right) - h\left(\frac{n_2}{p_i}\right) \geq 0. $$

We have $$h\left(\frac{n_1}{p_i}\right) - h\left(\frac{n_2}{p_i}\right) = h\left(\frac{m}{p_i}\right)\left(\frac{p_{\ell}^2+ p_\ell + 1}{p_\ell^2} -  \frac{p_\ell +1}{p_\ell}\frac{p_k+1}{p_k}\right). $$
Condition 1 guarantees that the above expression is positive.
 
\end{proof}

We have as an immediate consequence:

\begin{lemma} Assume $n$ is a record setter for $v$. If $17|n$ then $3^2|n$. If $37|n$, then $5^2|n$. If $67|n$, then $7^2|n$.\label{when small squares divide n}
\end{lemma}

\begin{lemma} If $p_i$ is the $i$th prime number for $i = 1 ,2 \cdots$ and $k \geq 29$, then
$$ \frac{2}{1} \frac{3}{2} \frac{5}{4}\cdots  \frac{p_k}{p_k-1} \leq  e^\gamma (\log p_k)\left(1+\frac{1}{1.2\log^2 p_k}\right),$$ where $\gamma$ is Euler's constant. \label{Modified RS Merten upper bound}
\end{lemma}
\begin{proof}
It follows from Theorem 8 of \cite{RosserSchoenfeld1962} that if $p_k \geq 293$, then
$$ \frac{2}{1} \frac{3}{2} \frac{5}{4}\cdots  \frac{p_k}{p_k-1} \leq  e^\gamma \left((\log p_k)\left(1+\frac{1}{2\log^2 p_k}\right)\right).$$
We then verify that the inequality holds for the primes between 29 and 293. 
\end{proof}

\begin{lemma} If $n$ has at least twenty distinct prime divisors, and we set $k=\omega(n)$, and $t=\log k$ then we have \begin{equation}\prod_{p|n} \frac{p}{p-1} \leq e^\gamma (\log \left(k \left(t + \log t -0.5 \right)\right)\left(1 + \frac{1}{1.2 \left(\log\left(k \left(t + \log t -0.5 \right) \right)\right)^2}\right)   \end{equation}
\label{Mert product upper bound}
\end{lemma}
\begin{proof} Assume that $n$ has $k=\omega(n)$ distinct prime factors. Since $\frac{x}{x-1}$ is a decreasing function for $x>1$, we may have $\prod_{p|n} \frac{p}{p-1}  \leq \frac{2}{1} \frac{3}{2} \frac{5}{4}\cdots  \frac{p_k}{p_k-1}$, where $p_k$ is the $k$th prime.

We have from Lemma \ref{Modified RS Merten upper bound} that  \begin{equation}\prod_{p|n} \frac{p}{p-1} \leq \label{upper bound on Mertens product} e^\gamma (\log p_k)\left(1+\frac{1}{1.2\log^2 p_k}\right). \end{equation}

Theorem 3 of Rosser and Schoenfeld \cite{RosserSchoenfeld1962} showed that for $k \geq 6$ the $k$th prime satisfies 

\begin{equation}p_k < k(\log  k + \log \log k -0.5).\label{RS upper bound on pk} \end{equation}

Combining Equation (\ref{upper bound on Mertens product}) and \ref{RS upper bound on pk} then gives the desired result.

\end{proof}

\begin{lemma} Let $p_i$ be the $i$th prime number. Then 
$\sum_{i=1}^k \frac{\log p_i}{p_i-1} \leq 1+ \log p_i $.
\label{RS derived estimate for sum of log p / p}
\end{lemma}
\begin{proof} This follows from estimates (2.8), (2.11) and (3.24) in \cite{RosserSchoenfeld1962}.
\end{proof}

\begin{lemma} If $n$ has at least $k$ distinct prime divisors for $k \geq 20$, then \begin{equation} \sum_{p|n} \frac{\log p}{p-1} \leq 1+ \log \left(k(\log  k + \log \log k -0.5)\right). 
\end{equation} \label{sum log p / p p divides n bound}
\end{lemma}
\begin{proof} This follows from combining \ref{RS upper bound on pk}  with Lemma \ref{RS derived estimate for sum of log p / p}.
\end{proof}

\begin{proposition} Any record setter for $v(n)$ has at most 29 distinct prime divisors. \label{Weak initial bound on number of distinct prime divisors of a record setter}
\end{proposition}
\begin{proof} Assume we have a record setter $n$ with at least $k \geq 19$ distinct prime divisors. Then by Lemma \ref{when small squares divide n}, $4|n$, $9|n$, $25|n$, and $49|n$. 

Thus $\tau(n) \geq 2^{k-2}(3^4)$ which is the same as 
$\log \tau(n) \geq (k-2)\log 2 + 4 \log 3$.

Thus we apply Lemma \ref{Lemma 4.3 from Weber}, Lemma \ref{Mert product upper bound}, and  Lemma \ref{sum log p / p p divides n bound} and set $\log k=t$, and $a = k(t + \log t -0.5)$ to obtain that

\begin{equation} v(n) <  \frac{\left(1+ \log a \right)e^\gamma (\log  a\left(1 + \frac{1}{1.2 \left(\log a \right)^2}\right)}{(k-2) \log 2 + 4 \log 3}.  \label{Long equation}
\end{equation}

Write the right-hand side of (\ref{Long equation}) as  $f(k)$. Note that the derivative of $f(k)$ is negative and so $f(k)$ is decreasing in $k$. Since we have that $$v(321253732800)=1.70595781\cdots,$$ if we have a record setter, we must have 
$f(k) > 1.7059578$. But this can only occur when $k< 54.4$, and since $k$ is an integer, we obtain that $k \leq 54$. Direct calculation then confirms that if $30 \leq k \leq 54 $, then $$\frac{z(n)}{(k-2)\log 2 +4 \log 3}< 1.7059578, $$
which completes the proof. 
\end{proof}

\begin{theorem} The largest record-setter for $v(n)$ is 321253732800.
\end{theorem}

\begin{proof} Our method of proof can be thought of as ``reverse bootstrapping.'' By Lemma \ref{Weak initial bound on number of distinct prime divisors of a record setter}, any record setter must have at most 29 distinct prime divisors. We have from Lemma \ref{Tightened Weber 4.3},

Assume for now that $n$ has exactly 29 distinct prime factors. We must have that $4|n$, $9|n$, $25|n$, and $49|n$. If we had $11^2|n$, then we would have $\tau(n) \geq 3^5 2^24$, and so $\log \tau(n) > 22.12859$, but then we would have $v(n) \leq \frac{37.083}{22.128} < 1.7$. Thus, 11 and other primes cannot be raised to more than the first power. If 3 is raised to a power higher than the second, then we must have 
$\tau(n) \geq (2^3)(3^3)2^25$. and so $\log \tau(n) \geq 22.7039$ , and so similarly we have then $v(n) \leq \frac{37.083}{22.7039}$. 

Thus, our only option here is that $n=2^23^25^27^2 (11)(13)(17) \cdots (109)$, which we can check. We can continue this way algorithmically at each stage identifying a finite list to check and verifying that they do not match. The code for this, and its results can be found at  \cite{Tcode}.

\end{proof}

\section{Lower bounds for $z(n)$}

So far, most of our results have focused on upper bounds for $z(n)$. Let us now turn to lower bounds for $z(n)$.

Assume that $n >1$ and let $p$ be the smallest prime divisor of $n$. Then we trivially have 
that $$z(n) = \sum_{d|n} \frac{\log d}{d} = \sum_{d|n, d>1} \frac{\log d}{d} \geq \frac{\log p}{d} = (\log p) \sum_{d|n,d>1} \frac{1}{d} = (\log p)h(n). $$

Here, again, the weighted AM-GM inequality can give us tighter bounds, albeit only marginally better.  As before, we will first see what we get when $n$ a perfect number, where the form is a little nicer. 
We choose as our choice of weights $\alpha_i = \frac{1}{d_i}$ with $1 \leq i \leq t$. But this time we also set $x_i= \frac{1}{d_i}$. Since $n$ is perfect, the sum of the weights is again 1. 

\begin{theorem} Let $n$ be a perfect number. Then,
$$\log \sum_{d|n, d>1} \frac{1}{d^2} \geq -  \sum_{d|n, d>1} \frac{\log d}{d} = - z(n).$$ \label{sum 1/d2 version}
\end{theorem}

We have as an immediate consequence the following corollary. 




We get a similar statement for pseudoperfect numbers.

Note that there is a massive gap here between our lower and upper bounds for $z(n)$ even in this case where $n$ is a perfect number.  In particular, it follows from results of Norton \cite{Norton}, that if $n$ is a perfect or abundant number with smallest prime factor $p$, then $n$ has at least roughly $\frac{p^2}{\log p}$ distinct prime factors, and thus we have roughly that  $2^\frac{p^2}{\log p} \leq \tau(n)$. 

As before, we get a similar inequality for all $n$. Setting $\alpha_i$ as running over $\frac{1}{d_i}$ where $d_i$ ranges over the divisors of $n$ which are greater than 1, and setting $x_i = \frac{1}{d_i\left(h(n)-1\right)}$ we obtain:

\begin{proposition} For all positive integers, $n>1$,  with smallest prime factor $p$, we have \label{z(n) general lower bound}
$$z(n) \geq \left(h(n)-1\right) \log \sum_{d|n, d>1} \frac{1}{d^2},$$ with equality if and only if $n$ is prime. We also have
$$z(n) \geq \left(h(n)-1\right) \log (2p).$$
\end{proposition}

It is not hard to improve Proposition \ref{z(n) general lower bound} slightly, using the Prime Number Theorem to get:

\begin{proposition} For any $\epsilon>0$, if $p$ is a sufficiently large prime, then for any $n$ with smallest prime factor $p$, $$z(n) \geq \left(h(n)-1\right) \log (2p (1-\epsilon) \log p).$$ \label{teeny tiny improvement to z lower}
\end{proposition}

Proposition \ref{teeny tiny improvement to z lower} can be made explicit with only a little work. Unfortunately, this is only a very tiny improvement since $(2p\log p)$ is itself inside a logarithm. In fact,
we cannot hope to improve these bounds more than a small amount without some additional information about $n$, since in Proposition \ref{aj has z not much more than log pj}, we exhibited an infinite sequence $c_j$ with smallest prime factor $p_j$ and  where $z(c_j) = O(\log p_j).$



Note that one cannot get a general lower bound for $z(n)$ in terms of $\tau(n)$ even if one fixes the smallest prime factor of $n$ to be $p$. In particular, fix $t$ to be a positive integer. Then for any $\epsilon >0$, there is an $n$ is with $\tau(n)=t+1$ and where 
$z(n) \leq(1+\epsilon)z(p^{q-1})$ where $q$ is the smallest prime factor of $t$. To see this, consider numbers $n$ of the form $p^{q-1}r^{\frac{t}q}$ where $r$ is a large prime.

Note also that it follows from Mertens' theorem that for any $\epsilon >0$ and positive constant $K$, there is an integer $n$ where $\frac{z(n)}{h(n)} >K$ and $h(n) < 1+\epsilon$. 

\section{Other weights}

Another weighing option for the weighted AM-GM inequality is to restrict our attention to perfect numbers and use all of the divisors, including $n$ itself, and then divide by $2n$.  We can use as weights $\alpha_i=\frac{d_i}{2n}$, and set $x_i =\frac{n}{d_i}$, where $i$ ranges from $0$ to $t$. 


Note that to have a chance at get an interesting result here we do need to that $n$ is a perfect number rather than just a pseudoperfect number, since otherwise $\frac{n}{d}$ will not be guaranteed to be a divisor in our divisor set.

From this choice of $x_i$ and $\alpha_i$ we obtain that if $n$ is a perfect number then, 

\begin{equation} \label{ai = d/2N, xi=N/d  } \log \frac{\tau(n)}{2} \geq \sum_{d|n}\frac{d \log \left(\frac{n}{d}\right)}{2n}. 
\end{equation}


%

If $n$ is perfect, a little work can show that Inequality (\ref{ai = d/2N, xi=N/d  }) is equivalent to Inequality (\ref{tau lower bound for perfect form}) by substituting $n/d$ for $d$. However, Inequality  (\ref{ai = d/2N, xi=N/d  }) turns out to be less interesting, because it is actually true for all $n$ for essentially trivial reasons since the function $f_n(x)=x \log \frac{n}{x}$ has its global maximum at $m=\frac{n}{e}$, where $f(m)=\frac{n}{e}$, so all the terms on the right hand side of Inequality (\ref{ai = d/2N, xi=N/d  }) are less than $\frac{1}{2}$.

Some other choices of weights are more interesting.  

Given the close connection between the Euler function $\phi(n)$ and $\sigma(n)$, applications involving $\phi(n)$ are a natural thing to investigate.  We may use the weighted AM-GM inequality with $\alpha_i = \frac{d}{2n}$ and $x_i = \phi(d)$.
This gives us:

\begin{proposition} If $n$ is a perfect number then $$\sum_{d|n} \frac{d}{2n} \log \phi(d) \leq \log \sum_{d|n} \frac{\phi(d)d}{2n}. $$  \label{A < B}  
\end{proposition}

We also note that by replacing $d$ with $\frac{n}{d}$, we obtain that $$\log \sum_{d|n} \frac{\phi(d)d}{2n} =  \log \sum_{d|n} \frac{\phi(\frac{n}{d})}{2d}.$$ Thus, one also has:

\begin{proposition} If $n$ is a perfect number then $$\sum_{d|n} \frac{d}{2n} \log \phi(d) \leq \log \sum_{d|n} \frac{\phi(\frac{n}{d})}{2d}. $$    
\end{proposition}

We wish to extend these notions to pseudoperfect numbers.   Given a positive integer $n$ and $S$, a subset of the divisors of $n$, we will define three functions as follows. We define

\begin{equation}A(S,n) = \sum_{d \in S } \frac{d}{2n} \log \phi(d), \label{A def}\end{equation}

\begin{equation} B(S,n) =  \log \sum_{d \in S} \frac{\phi(d)d}{2n}, \label{B def}
\end{equation}

and 

\begin{equation} C(S,n)=  \log \sum_{d \in S} \frac{\phi(\frac{n}{d})}{2d} \label{C def}.    
\end{equation}

We immediately have from using the same weighing as in Proposition \ref{A < B} that for any pseudoperfect $n$ and corresponding set $S_1$ we have $A(S_1,n) \leq B(S,n)$. However, it is not true that $C(S_1,n) = B(S_1,n)$ for any pseudoperfect $n$ and corresponding set $S_1$. The proof for perfect numbers used the replacement of $d$ with $\frac{n}{d}$ and for a general pseudoperfect number this cannot happen. It is possible, that $d \in S$ but $\frac{n}{d} \not \in S$. In these circumstances, $C(S_1,n)$ can not only be less than $A(S_1,n)$ it can even be negative. For example,  $C(\{13, 26,39, 78,\},78) =   -2.0053 \cdots$. 

One can reasonably object that the example of $n=78$ is extreme, and that this number is obviously not really a perfect number since the set does not include 1. And in fact, it is not hard to show that if $1 \in S_1$, then $C(S_1,n) > 0$. However, there are examples of pseudoperfect $n$ which have an $S_1$, where $1 \in S_1$ and yet $A(S_1,n) > C(S_1,n)$.  An example is $n=90$ with the set $S_1=\{1, 3, 5, 6, 30, 45, 90\}.$ 

We can show even more pathological behavior. Consider $n=272$ with $$S_1= \{1, 16,17, 34, 68, 136, 272\}.$$  This example is interesting because $272$ is a primitive non-deficient number, which is a much stronger property than being a primitive pseudoperfect. Of course, even primitive non-deficients are much more common than odd primitive non-deficients. However, this example has $A(S_1,n) > C(S_1,n)$.

 \begin{question} Is there an odd primitive non-deficient $n$ which is pseudoperfect with set $S_1$, with $1 \in S_1$, and where $A(S_1,n) > C(S_1,n)?$     
 \end{question}

 Note that $n=1575$ almost provides an answer to the above question.  In particular, $n=1575$ is an odd primitive non-deficient, pseudoperfect number with pseudoperfect set $$S_1=
 \{225, 35, 5, 7, 9, 15, 21, 63, 75, 105, 175, 315, 525,    1575 \} $$ where $A(S_1,n) >C(S_1,n)$. However, $1$ is not an element of $S_1$.

However, note that if a pseudoperfect number $n$ does satisfy that $d$ is in its sum set if and only if $\frac{n}{d}$ is in its sum set, then  for that choice of $S_1$ we will have $B(S_1,n)= C(S_1,n)$ and hence they will  satisfy that $A(S_1,n) < C(S_1,n)$. Motivated by this observation, we will say that a number $n$ is \emph{strongly pseudoperfect} if it is a pseudoperfect number with a divisor sum set $S_1$ where $d \in S_1$ if and only if $n/d \in S_1$.  The strongly pseudoperfect numbers are given by A334405 in the OEIS\cite{OEIS A334405}, and the first few are 
$$6, 28, 36, 60, 84, 90, 120, 156, 210, 216, 240, 252 \cdots .$$

A prior paper \cite{AMPSZ} work by Aryan, Madhavani, Parikh,  Slattery, and the second author of this paper classified strongly pseudoperfect $n$ where $n=2^kp$ for some odd prime $p$ and where exactly a single divisor pair is omitted from the set of divisors. 

Strongly pseudoperfect numbers can be thought of as a generalization of Descartes' ``spoof perfect number.'' Descartes noticed the number $$D=(3^2)(7^2)(11^2)(13^2)(22021).$$  He saw that $D$ looks almost like it is an odd perfect number, in that if one naively applies the sum of divisors formula to $D$, one has 

\begin{equation}\sigma(D)=(3^2+3+1)(7^2+7+1)(13^2+13+1)(22021+1) = 2(3^2)(7^2)(11^2)(13^2)(22021).\label{Descartes spoof} \end{equation}

A reader who has not previously seen Equation (\ref{Descartes spoof}) should be immediately suspicious because all the factors on the left-hand side are much smaller than the supposed prime 22021. (The exception is of course $2202+1$, but $22021$ and $22021+1$ must be relatively prime.) And in fact, 22021 is not prime. But if we pretend that it is prime, then we would have an odd perfect number.  Spoof perfect numbers are an obstruction to proving the non-existence of odd perfect numbers because many non-trivial statements one can prove about odd perfect numbers will also have corresponding, nearly identical statements about spoofs. Since spoof perfect numbers exist, any collection of such statements will then be insufficient to prove the non-existence of odd perfect numbers.

Any spoof perfect number in the sense of Descartes will give rise to a strongly pseudoperfect number if one expands out the factorization, and using the factorization implied by pretending that the corresponding factors are prime. Some authors have referred to the set of all spoof perfect numbers and actual perfect numbers as ``freestyle perfect numbers.''  What precisely is meant here has some degree of ambiguity. For example, when some $p$ is wrongly assumed to be prime, do we wish to only do so if none of factors of $p$ appear as primes themselves? See OEIS A174292 for some discussion. 

In recent years, spoof perfect numbers have been generalized and made more rigorous. John Voight introduced the example $V=(3^4)(7^2)(11^2)(19^2)(-127)$, where if we ignore the somewhat glaring issue that $-127$ is a negative number when we apply the usual formula for $\sigma$ we obtain that
$$(3^4+3^3+3^2+3+1)(7^2+7+1)(11^2+11+1)(19^2+19+1)(-127+1) = 2V.$$

Another way we can allow for spoofs, introduced in  \cite{BYU}, is where we are allowed to fail to notice that we have repeated a prime factor. For example, consider the number $$P=(3^2)(7^2)(7^2)(13^1)(-19)^2.$$

If one ignores that $7^2$ and $7^2$ are not powers of distinct primes, and ignores that $(-19)$ is negative (as with Voight's example), then one has that
$$(3^2+3+1)(7^2+7+1)(7^2+7+1)(13+1)((-19)^2 + (-19)+1)= 2P.$$ 

Continuing with this idea, \cite{BYU} gave a more rigorous expansion of these spoofs and showed that there are infinitely many and of a wide variety of forms. 

In any event, 156 is strongly pseudoperfect since the sum of all its divisors, excepting 2 and 78, adds up to 2(156), and 156 is not a freestyle perfect or spoof for any of these definitions. Thus,  being a strongly pseudoperfect is a broader class of numbers than spoofs.

The more general notions of spoof perfect numbers due to Voight and Nielsen can give rise to a corresponding and more general notion of strongly pseudoperfect numbers, where one is allowed to have negative divisors in $S_1$, and where $S_1$ rather than be a set, becomes a multiset. However, there is some room for hope that statements proven using the sort of techniques in this paper may still be of some interested for the perfect number problem. In particular, the weighted AM-GM inequality only applies to positive numbers, and Zaremba's function does not make sense for negative numbers. Although there are inequalities similar to the weighted AM-GM inequality which allow zero, they generally lose something. See for example, the Steffensen inequality \cite{Steffensen}. We will not investigate this broader notion here.

If one follows one Descartes's lead and restricts one's ``freestyle perfect numbers'' to only being allowed to pretend that a number is prime if that number is genuinely relatively prime to all the other primes (real or pretend) in the factorization, and do not allow negative primes, then every spoof perfect number gives a strongly pseudoperfect number. 


It is plausible that tightened versions of the weighted AM-GM inequality or related bounds which use that we are taking sums over genuinely distinct integers can at least avoid the obstruction posed by the most general form of spoofs. For example, our argument giving a lower bound for $z(n)$ used that the sum of $\frac{1}{d^2}$ was ranging over distinct positive integers $d$, and our proof would not apply if repetition of divisors is allowed.\\

Unfortunately, all of these results apply to all Cartesian spoofs. Worse, when one has applied Mertens' theorem or similar density statements, it seems unlikely that they will be highly informative, since prime density results cannot see if one has slipped in a single non-prime or even a few non-primes if they are sparse enough. \\

With that pessimistic digression out of the way, let us at least see what we can prove about strongly pseudoperfect numbers.

\begin{proposition}    
If $m \geq 2,$ then $n= 2^{m-1}(2^m-1)$ then $n$ is strongly pseudoperfect. \label{Almost Euclid}
\end{proposition}
\begin{proof} The proof is essentially the same as the proof that if $2^p-1$ is prime, then $2^{p-1}(2^p-1)$ is perfect, but we get to ignore any unwanted contribution from divisors of $2^n-1$. 
\end{proof}

Note that all the numbers constructed by Proposition \ref{Almost Euclid} are even. We suspect that there are infinitely many odd numbers which are strongly pseudoperfect but do not see an obvious way to prove it. 

If $n$ is pseudoperfect, then so is $mn$ for any $m$, since one take each element in $S_01$  and multiply it by $m$. However, this is not the case for strongly pseudoperfect numbers. The next two results show different ways this can fail.

\begin{proposition} If $n$ is strongly pseudoperfect, then $n$ is not $3$ (mod 4) and $n$ is not $2$ (mod 3). \label{strong pseudos satisfy Touchard}
\end{proposition}
\begin{proof} We will prove that $n$ is not 2 (mod 3). The proof for $3$ (mod 4) is similar. Suppose that $n$ is strongly pseudoperfect with divisor set $S$, and that $n$ is 2 mod 3. Since $n$ cannot be a perfect square, we may pair up the divisors $d$ in $S$, as $(d, \frac{n}{d})$. Note, if $d \equiv 1$ (mod 3), then $\frac{n}{d} \equiv 2$ (mod 3), and that if $d \equiv 2$ (mod 3), then $\frac{n}{d} \equiv 1$ (mod 3). No divisor can be $0$ (mod 3), because $n$ is not divisible by $3$.  Thus, for each such pair we have $d+ \frac{n}{d} \equiv 0$ (mod 3), and so the total sum of the divisors is $0$ (mod 3). This is a contradiction since the sum of all the divisors is $2n$ which is $1$ (mod 3).
\end{proof}

Proposition \ref{strong pseudos satisfy Touchard} shows, incidentally, that if $n$ is either an odd strongly pseudoperfect number, or is a strongly pseudoperfect number which is not divisible by 3, then a positive proportion of its multiples will fail to be strongly pseudoperfect. We suspect that the modulus restrictions in Proposition \ref{strong pseudos satisfy Touchard} are the only modulus restrictions on strongly pseudoperfect numbers. 

If the Artin primitive root conjecture is true, (in particular if $2$ is a primitive root for a positive proportion of primes), then Proposition \ref{Almost Euclid} would imply that a positive proportion of primes have a strongly pseudoperfect number in at least half their congruence classes.\\
One may be tempted to wonder if $n$ is a strongly pseudoperfect number which is divisible by 6, will all multiples of $n$  be strongly pseudoperfect? The next result shows that this is not the case.

\begin{proposition} Let  $n$ be a positive integer and let $P$ be a prime with $P> \sigma(n)$. Then $nP$ is not strongly pseudoperfect. \label{big class of non- strongly pseudoperfect}
\end{proposition}
\begin{proof} Note that we have $$\sum_{d|nP, (d,P)= 1}d = \sigma(n).$$ Thus, if $nP$ is strongly pseudoperfect with set $S_1$, then at least one element in $S_1$ must be divisible by $p$, and $S_1$ must contain at least one divisor of $n$. We can write 
$$\sum_{d \in S_1} d = \sum_{d \in S_1, (d,P)=1} d + \sum_{d \in S_1, P|d}d $$ where both sums on the right-hand side are positive. But the second sum is divisible by $P$ whereas the first is not, so the sum of all elements in $S_1$ cannot be divisible by $P$, hence $\sum_{d \in S_1} d \neq nP.$   
\end{proof}



An open question of Erd\H{o}s \cite{Benkoski and Erdos} is whether there is a constant $C$ such that if $h(n) > C$, then $n$ must be pseudoperfect. Proposition \ref{big class of non- strongly pseudoperfect} shows that this is not the case for strongly pseudoperfect numbers. The pseudoperfect numbers have positive natural density, but it seems likely that the strongly pseudoperfect numbers have natural density zero.

Aside from the Euclid-like strongly pseudoperfect numbers, explicit families of strongly pseudoperfect numbers seem hard to come by. However, we note that $60$, $120$, $240$, $480$, $960$, $1920$, $3840$, $7680$, $15360$ are all strongly pseudoperfect.  

\begin{question} For any $n \geq 0$, is $60(2^n)$ strongly pseudoperfect?
\end{question} 

There may be other similar families. However, any $m$ where $m(2^n)$ is strongly pseudoperfect for all $n$ must have $3|m$. To see this, note that  if $m \not \equiv 0$ (mod 3) then either $m \equiv 2$ (mod 3), or $2m \equiv 2$ (mod 3). Thus at least one of this pair of numbers must fail to be strongly pseudoperfect. So hope for any similar family where we replace $60$ with another $m$ must at least have $m \equiv 0$ (mod 3). This condition is necessary but it is not sufficient: 90 is strongly pseudoperfect but 180 is not.\\

Note also that our definition of strongly pseudoperfect does not insist that $n \in S_1$. Any genuinely perfect number will have $1$ and $n$ in its corresponding set of divisors which sums to $2n$. Define a number to be extremely strongly pseudoperfect if it is strongly pseudoperfect for some set $S_1$ with $n \in S_1$. 

\begin{question} Is every strongly pseudoperfect number extremely strongly pseudoperfect?
\end{question}

Pollack and Shevelev \cite{PollackShev} defined a class of pseudoperfect numbers they called ``near-perfect numbers.''  A number $n$ is said to be near-perfect number when the sum of all $n$'s proper divisors except a single divisor $k$ adds up to the number. Further investigation of near perfects includes \cite{RenChen}. An example is 12, where the proper divisors of $12$ are $1,2,3,4,$ and $6$, and one may omit $k=4$ to have $1+2+3+6=12$. One might expect that no strongly pseudoperfect number can be a near perfect number; here one removes only a single divisor and a strongly pseudoperfect number must lose at least two. \\

However, there are two problems with this line of reasoning. First, it may be that one has a near perfect number $n=k^2$ where the single divisor omitted is $k$, in which case the same set would give a strongly pseudoperfect sum for $n$. Second, it may be that the omitted divisor is itself the sum of a pair of divisors, and so we can take out that pair of divisors instead of $k$. Notice that for the near perfect number 12, one could instead of omitting 4, omit the pair of divisors $1$ and $3$ to get a pseudoperfect sum for 12. However, since $(1)(3) \neq 12$, this does not give rise to a strongly pseudoperfect  representation.

Pomerance and Shevelev constructed what is likely an infinite family of perfect square near-perfect numbers, since they showed that $2^{p-1}(2^p-1)^2$ is near perfect, whenever $(2^p-1)$ is prime. However, the omitted divisor $k$ is $2^p-1$, not the square root. This prompts three questions.

\begin{question} Is there any near-perfect number $n$ which is a perfect square, the omitted divisor $k$ is in fact the square root of $n$? \label{near perfect with perfect square as omitted} Equivalently, does there exist an $m$ such that $\sigma(m^2)=2m^2 +m$?
\end{question}

\begin{question} Are there infinitely many strongly pseudoperfect perfect squares? \label{any strong pseudoperfect perfect square}
\end{question}

\begin{question} Is there any number which is both near-perfect and strongly pseudoperfect? \label{any near perfect and strongly pseudoperfect}
\end{question}

Note that affirmative answers to Question \ref{any strong pseudoperfect perfect square} and Question \ref{near perfect with perfect square as omitted} would imply an affirmative answer to Question  and \ref{any near perfect and strongly pseudoperfect}. 

Note while this paper was undergoing revisions, Question \ref{any strong pseudoperfect perfect square}, was found to have an affirmative answer by a student of the second author, in a publication to appear \cite{WZ}. 

It is known that pseudoperfect numbers have positive natural density  and that perfect numbers have zero density.\cite{Benkoski and Erdos}\cite{Hornfeck and Wirsing} Let $SP$ be the set of strongly pseudoperfect numbers.  

\begin{question} Does $SP$ have positive natural density? \label{strongly pseudoperfect natural density}
\end{question}

Question \ref{strongly pseudoperfect natural density} is connected to another obstruction to proving that no odd perfect numbers exist. In \cite{Zelinskybig}, the second author introduced the following idea. Given a set of natural numbers $S$, we will write $S(x)$ as the number of elements in $S$ which are at most $x$.  Euler proved that if $n$ is an odd  perfect number, then $n=q^sm^2$ for some prime $q$, where $q \equiv s \equiv 1$ (mod 4), and where $(q,m)=1$. Let $E$ be the set of all numbers of that form. Given a statement $S$ we will write $E_S$ to be the set of elements of E which satisfy $S$. The second author defined a property $S$ to be a weak property if
$$\lim_{x \rightarrow \infty} \frac{E_S(x)}{E(x)} = 1.$$ We observe that if one has a finite list of properties which an odd perfect number must satisfy, and all of them are weak properties, then they cannot suffice to prove there are no odd perfect numbers. That is, almost all elements in $E$ will not have been ruled out, as the intersection of all the relevant sets will still have positive density. Motivated by this, the second author said that $S$ is a strong property if $$\lim_{x \rightarrow \infty} \frac{E_S(x)}{E(x)} = 0.$$

In that context, it seems natural to ask which properties proven here are weak or strong properties. For example, simply being non-deficient is not a strong property, but it is not weak either. In particular, let $D_0$ be the statement ``N is not deficient.'' Then it is not too hard to show that.

$$0< \lim_{x \rightarrow \infty} \frac{E_{D_0}(x)}{E(x)} < 1.$$

We then ask:

\begin{question} Is being a strongly pseudoperfect number a strong restriction, a weak restriction, or neither?
\end{question}

Similarly, let $D_1$ be the statement ``n is pseudoperfect with some $S_1$ where $A(S_1,n) \leq C(S_1,n)$.''

\begin{question} Is $D_1$ a strong property, a weak property, or neither?
\end{question}

Another variant of perfect numbers is that of a \emph{primary pseudoperfect number}.\cite{Sondow and MacMillan}  A number $n$ is said to be a primary pseudoperfect number if $n$ satisfies the equation
$$\frac{1}{n}+\sum_{p|n} \frac{1}{p} =1,$$ where the sum is taken over all prime divisors of $n$. For example, $n=42$ is primary pseudoperfect since one has $$\frac{1}{42} + \frac{1}{2} + \frac{1}{3} + \frac{1}{7} =1. $$ It is not hard to see that if $n$ is primary pseudoperfect and $n$ is not $1,$ $2,$ or $6$, then $n$ must be abundant. Primary pseudoperfect numbers are very rare. We suspect the following.

\begin{conjecture} \label{no primary pseud is strongly pseud}
6 is the only primary pseudoperfect numbers which is also strongly pseudoperfect. \end{conjecture}

A number is called a \emph{weak primary pseudoperfect numbers} if it satisfies the weaker condition that $$\frac{1}{n}+\sum_{p|n} \frac{1}{p} \in \mathbb{N}.$$ Note that now the sum is allowed to be any integer, not just 1. Other than $n=1$, there is no known example of a weak primary pseudoperfect number that is not in fact a primary pseudoperfect number. It may seem safe then to also make what seems on its face to be only a slightly stronger conjecture.

\begin{conjecture} \label{no weak primary pseud is strongly pseud}
6 is the only weak primary pseudoperfect numbers which is also strongly pseudoperfect. \end{conjecture}

At first glance, Conjecture  \ref{no weak primary pseud is strongly pseud} to be only a hairsbreadth stronger than Conjecture \ref{no primary pseud is strongly pseud}, but this may not be the case. It is plausible that there are many weak primary pseudoperfect numbers that are not primary pseudoperfect as one gets to large numbers. In any event, we suspect that that Conjecture \ref{no weak primary pseud is strongly pseud}, if true, will be substantially more difficult to prove than Conjecture  \ref{no primary pseud is strongly pseud}.

Closely related to the weak primary pseudoperfect numbers are the \emph{Giuga numbers}, which satisfy the similar equation $$\frac{-1}{n}+\sum_{p|n} \frac{1}{p} \in \mathbb{N}.$$ The first few Giuga numbers are 30, 858, 1722, 66198 and are OEIS sequence A007850.\\

Giuga numbers were extensively studied in  \cite{Borwein}. The original motivation for studying Giuga numbers was the Agoh–Giuga conjecture which states that if $n >1$, then $n$ is prime if and only if $nB_{n-1} \equiv -1$ (mod $n$). Here $B_k$ is the $k$th Bernoulli number.  The Agoh–Giuga conjecture is equivalent to the statement that if $n > 1$, then $n$ is prime if and only if $$1^{n-1}+2^{n-1}+ \cdots +(n-1)^{n-1} \equiv -1 \pmod n.$$ Any counterexample to the Agoh-Giuga conjecture must be a Giuga number; this fact was the original motivation for studying Giuga numbers. We suspect that there are no strongly pseudoerfect Giuga numbers. 

Recently, Grau, Oller-Marc\'en, and Sadornil\cite{GOS}, generalized the notion of weak primary pseudoperfect numbers and of Giuga numbers to what they called a \emph{$\mu$-Sondow number}, which is a number $n$ along with an integer $\mu$, such that 
\begin{equation}\frac{\mu}{n}+\sum_{p|n} \frac{1}{p} \in \mathbb{N}. \label{mu Sondow form 1}\end{equation}

Equation (\ref{mu Sondow form 1}) is equivalent to the statement that 

\begin{equation}\sum_{p|n} \frac{n}{p} \equiv -\mu \pmod n, \label{mu Sondow form 2}\end{equation}

Therefore every $n$ is a $\mu$-Sondow number for some choice of $\mu$.

Grau, Oller-Marc\'en, and Sadornil conjectured that for any  $\mu$, there is a $\mu$-Sondow number greater than $|\mu|$ but noted that for any fixed $\mu$,  $\mu$-Sondow numbers are rare.  In this context, we have the following conjecture:

\begin{conjecture} \label{Finite strongly pseudoperfect for any mu}
For any $\mu$, there are only finitely many $\mu$-Sondow numbers which are strongly  pseudoperfect. \end{conjecture}

Note also that 11025 is a curious number for being odd, strongly pseudoperfect and for being primitive non-deficient (that is, it is not deficient and has no non-deficient proper divisors). Are there infinitely many such numbers?










For much of this paper, we have assumed that $n$ is perfect or pseudoperfect number and used that fact to construct a choice of weights. However, one of the basic facts about $\phi(n)$ is that for any $n$, it satisfies $\sum_{d|n} \phi(d) =n$. In that context, another choice of weights is $\alpha_i = \frac{\phi(d)}{n}$ where $d$ ranges over the divisors of $n$. If one uses this choice with $x_i=d$, then the weighted AM-GM inequality yields the following:

\begin{proposition} For any positive integer $n$, $$\sum_{d|n} \frac{\phi(d) \log d}{n} \leq \log \sum_{d|n} \frac{\phi(d)d}{n}.$$
\end{proposition}

Note that this means that when $n$ is perfect, we have two natural sets of numbers which both sum 1, namely $\frac{\phi(d)}{n}$ and $\frac{1}{2d}$. One can think of these as two different probability distributions on set of divisors, and can use ideas from probability to derive other inequalities. The Kullback–Leibler divergence is one of the more obvious choices here, and a sequel paper to this one will discuss the consequences of applying it in this and related contexts.

\section{Tightened versions of the AM-GM Inequality}

Most of the bounds here rely  on appropriately chosen weights for the weighted AM-GM inequality. This raises a question: There are many tightened versions of the inequality which specifically apply when we know more about the set of numbers and choice of weights. Can we take advantage of them here? The answer appears to be ``yes, but not by much.'' However, it is plausible that further tightened versions that use more information about the sets in question,  would work. But the ones in the literature currently are not well-optimized for our purposes. We will discuss some of these, and present a modified versions specifically using information from the situations of interest. However,  none of these will add more than a small improvement to our inequalities from earlier without some further insight.


We list one potential exception here. Aldaz \cite{Aldaz} has the following version of the weighted AM-GM inequality.

\begin{equation} \label{Aldaz AM-GM}
\prod_{i=1}^n x_i^{\alpha_i} \leq \sum_{i=1}^n{\alpha_i x_i} - \sum_{i=1}^n \alpha_i\left(x_i^{\frac{1}{2}} - \sum_{j=1}^n \alpha_j x_j^{\frac{1}{2}} \right)^2.
\end{equation}
When applying Equation (\ref{Aldaz AM-GM}) with our earlier weighing to with $\alpha_i = \frac{1}{d_i}$ where we range over the non-1 divisors of a perfect number $n$, we obtain

\begin{equation} \prod_{d|n}d^{\frac{1}{d}} < \tau(n)-1  - \sum_{d|n, d>1} \frac{1}{d}\left(d^{\frac{1}{2}} - \sigma_{\frac{-1}{2}}(n) +1 \right)^2. 
\end{equation}


It seems plausible that a non-trivial estimate of the last term could lead to an interesting improvement in the general inequality for the perfect case. A somewhat sharper inequality also due to Aldaz is in \cite{Aldaz2}.


One approach of interest is to say that the weighted AM-GM inequality is just essentially Jensen's inequality with $f(x) = \log x$. We can tweak our choice of $f(x)$, using that for our purposes, $f(x)$ needs only obey a concavity or convexity requirement in a specified interval. Since $n$ is generally much larger than its smallest prime factor $p$, if one is going to use that all non-1 divisors are in the range $[p,n]$, using the lower bound on the domain is more helpful than using the upper bound.

We will therefore choose $f(x) = \log x + \frac{A}{2}\frac{1}{x}$ which is chosen so that it is concave in the interval $[A,\infty]$.

\begin{lemma} Let $A$ be a positive real number, and let $f(x)= \log x + \frac{A}{2}\frac{1}{x}$. Then $f''(x) \leq 0$ for all $x$ in the interval $[A, \infty).$
\end{lemma}
\begin{proof} This follows immediately from the fact that $f''(x) = \frac{-1+ \frac{A}{x}}{x^2}$.
\end{proof}

We may then apply Jensen's inequality to $f(x)$ obtain the following.

\begin{lemma} Assume that $\alpha_1, \alpha_2 \cdots \alpha_k$ are positive real numbers where $\alpha_1 + \alpha_2 \cdots + \alpha_k =1$. Assume that $x_1, x_2 \cdots x_n$ are real numbers in the interval $[A, \infty]$. Then
\begin{equation} \label{GM-AM where we know the smallest value}
\log(\sum_{i=1}^k \alpha_i x_i) + \frac{A}{2}\left(\sum_{i=1}^k \alpha_i x_i\right)^{-1} \geq \sum_{i=1}^k \alpha_i \log x_i + \frac{A}{2}\sum_{i=1}^k  \frac{\alpha_i}{x_i}.
\end{equation}
\end{lemma}

From which we can obtain the following using our earlier weighing choice of $\alpha_i = \frac{1}{d_i}$, and $x_i = d_i$ to obtain:

\begin{proposition} If $n$ is an odd perfect number with smallest prime factor $p$, then 
$\log(\tau(n)-1) + \frac{p}{2(\tau(n)-1)} \geq z(n)  + \frac{1}{2p}$.
\end{proposition}


There is also some hope of using existing versions of the weighted AM-GM inequality which are already in the literature.

One variant of the weighted AM-GM inequality in the literature is due to Sababheh \cite{Sababheh} and is as follows. 

\begin{lemma} Assume that $\alpha_1, \alpha_2, \cdots \alpha_n$ are positive numbers with $\alpha_1 + \alpha_2 \cdots +\alpha_n=1$. Set $\alpha_{\mathrm{min}}$ to be the minimum of the $\alpha_i$, and $\alpha_{\mathrm{max}}$ to be the maximum of the $\alpha_i$. Let $x_1, x_2 \cdots x_n$ be positive numbers, and let $A$ be the unweighted arithmetic mean of the $x_i$, and $G$ be the unweighted geometric mean of the $x_i$ Then we have  \label{Sababheh's First Lemma}
\begin{equation} \label{Sababhe's first inequality}
\left(\frac{A}{G}\right)^{n\alpha_{\mathrm{min}}} \prod_{i=1}^n x_i^{\alpha_i} \leq \sum_{i=1}^n \alpha_i x_i \leq \left(\frac{A}{G}\right)^{n\alpha_{\mathrm{max}}} \prod_{i=1}^n x_i^{\alpha_i} .
\end{equation}
\end{lemma}

Using the same method as before, We obtain from Lemma \ref{Sababheh's First Lemma} the following.

\begin{proposition} If $n$ is a strongly pseudoperfect number with divisor set $S_1$, then $$\frac{|S_1|\log n}{3n} +\sum_{d \in S_1}  \frac{\log d}{d} \leq\log (|S_1|-1).$$ \label{strongly pseudoperfect using Sababheh}
\end{proposition}
\begin{proof}  We apply Lemma \ref{Sababheh's First Lemma} with $\alpha_i=1/(2d_i)$ running over $d_i$ and $x_i =d_i$ as before, and the proof is essentially identical to our earlier result, except we also use that the geometric mean of the elements in $S_1$ must be $\sqrt{N}$. Note that the proof that the geometric mean of the divisors of a number is equal to the square root was essentially using fact that divisors for non-squares can be paired up via $(d,\frac{n}{d})$. Thus, the same proof still goes through for the geometric mean of the $S_1$ set for a strongly pseudoperfect number, which completes the proof.
\end{proof}

In the specific case of perfect numbers Proposition \ref{strongly pseudoperfect using Sababheh} becomes:

\begin{corollary} If $n$ is a perfect number then $$\frac{\tau(n)\log n}{3n}  + z(n) \leq \log (\tau(n)-1).$$
\end{corollary}

But this is a very weak improvement. In particular,  $\frac{\tau(n)\log n}{3n}$ goes to zero as $n$ goes to infinity.

Another inequality which may be of interest is the weighted Ky Fan's inequality. This  stronger version of the weighted AM-GM inequality states that if $x_1, x_2 \cdots x_n$ are all non-negative numbers at most $\frac{1}{2}$, and $\alpha_1, \alpha_2 \cdots \alpha_n$ are non-negative real numbers satisfying our usual condition that $\alpha_1 + \alpha_2 \cdots \alpha_n=1$, then 

\begin{equation} \frac{ \prod_{i=1}^n x_i^{\alpha_i} }{ \prod_{i=1}^n (1-x_i)^{\alpha_i} }  \leq \frac{ \sum_{i=1}^n \alpha_i x_i }{ \sum_{i=1}^n \alpha_i (1-x_i) } \label{Weighted Ky Fan}.
\end{equation}

However, we have been unable to obtain a non-trivial result using the Ky Fan inequality which is substantially tighter than any of our results here. 
Similarly, a far reaching generalization of the Ky Fan inequality is Levinson's inequality. 

\begin{theorem}(Levinson)\cite{Levinson}
Let $a$ be a positive constant, and let $f(x)$ be a function whose third derivative is positive in $(0,2a)$. Assume further that 
$\alpha_1, \alpha_2 \cdots \alpha_n$ are positive real numbers with $\alpha_1 + \alpha_2 \cdots +\alpha_n =1$, and that $x_1, x_2 \cdots x_n$ are in the interval $(0,a)$. Then

\begin{equation}\label{Levinson inequality}
\sum_{i=1}^n \alpha_i f(x_i) -f\left(\sum_{i=1}^n\alpha_ix_i\right) \leq \sum_{i=1}^n\alpha_i f(2a-x_i)-f\left(\sum_{i=1}^n\alpha_i(2a-x_i)\right).
\end{equation}
\end{theorem}

 The conditions for Levinson's inequality have been further generalized in \cite{LawSeg}. 

In the cases of both Ky Fan's inequality and Levinson's inequality, there is a large amount of flexibility in the choices involved, and we suspect that it should be possible to use these inequalities to incorporate more number theoretic information into inequalities to get tighter results.  One potential goal for future is under the assumption that $n$ is a perfect to get a lower bound for $z(n)$ in terms of the number of distinct prime divisors of $n$. 



\appendix

\section{Algorithms and list of record setters}

In order to find all record-setters for $z(n)$ and $v(n)$ up to a maximum n, we generated all waterfall numbers up to n and then iterated through them in order, computing z and v. For this task, we found it most effective to treat waterfall numbers as a product of primorial numbers. Any product of primorials is a waterfall number and every waterfall number can be represented as a product of primorials, so algorithmically a waterfall number can be represented by a list of exponents of the successive primorials. For example, the list $[3, 0, 1]$ would represent $(2^3)(6^0)(30^1)=240$ The empty list represents the empty product, or 1.

The code starts with an empty list and uses a recursive scheme to construct the waterfall numbers from 0 up to and possibly including a max-$n$. Each call iteratively increments the last exponent in the list and outputs the resulting products until the product would exceed max-n. For each pass through the loop, the code also makes a recursive call that adds a 0 to the end of the list (until it reaches a primorial larger than max-$n$). In this way, the code generates every possible list of primorial exponents that represents a waterfall number that is within range. For details, see  \cite{Tcode}.

     \begin{center}
\begin{tabular}{ |c|c|c|c|c| } 
 \hline
$n$ & $z(n)$ & $\tau(n)$  & $v(n)=\frac{z(n)}{\log \tau(n)}$ & Record type  \\ 
1 & 0 & 1 & - & z\\
2 & 0.34657359028 & 2 & 0.5 & Both \\
3 & 0.36620409622 & 2 &  0.52832083357 & Both \\
4 & 0.6931471805599453 & 3 & 0.6309297535714574	& Both \\
6 &	1.0114042647073518 & 4 & 0.7295739585136225	& Both \\
12 &	1.5650534091363246 & 6 & 0.8734729387592397 & Both\\
24 & 1.9574025114441351	& 8	&  0.9413116320946855 & Both \\
36 &	 2.0693078747916944	&  9 & 0.9417825998016082 & Both \\
48	&  2.2113393276447026 &	 10 & 0.9603724676117413 & Both\\
60	&   2.6291351167661694 & 12 &  1.0580418049066245 & Both \\
120	&  3.1536019699500124 & 16	& 1.1374214807461371 & Both \\
180	&  3.296829727702829 &  18	& 1.140624806721235	& Both \\
240	& 3.485150114597928	&  20 &  1.1633716889070718	& Both \\
360	& 3.897045010404812	& 24 &  1.2262363126738316	& Both \\
720	& 4.272243596316464	& 30 &  1.2560998721657932 & Both \\	
840	& 4.438078029538006	&  32 & 1.2805586328585488	 & Both \\
1260 & 4.611033134346248 &  36 & 1.286733295829273 & Both \\
1680 &	4.844788911264986 & 40 & 1.3133497506571945 & Both \\	
2520 &  5.3572240096668935 &  48 &  1.3838661424637553 & Both \\
5040 &	5.816137669682257 &  60	& 1.4205296064593296 & Both \\
10080 &	6.08850361086746 &	 72 & 1.4236565215589612 & Both \\	
15120 & 6.248037988991858 & 80 & 1.425832955818364	& Both \\
25200 &	6.387503013476975 &	 90	&1.4195051527608638	& z \\
27720 &	6.653923297464633 &	 96 &  1.4578036158378216 & Both \\
50400 &	 6.675036433693168 &	 108 & 1.4256406131942176 & z\\
55440 &	7.181545678560509 & 120 & 1.5000643477637126 & Both\\	
110880 & 7.492166808574835 & 144 & 1.5075348623687792 &Both \\	
166320 &7.674162694562543 &160 & 1.5120984963170123 & Both\\	
277200 & 7.832742267027409 &180 & 1.5083395629764706 & z \\ 	
332640 & 7.99772386924453& 192 & 1.5212041672534768	& Both\\ 
554400 & 8.160359559426864 & 216& 1.518127795531729 & z\\
665280 & 8.183509553747768 & 224 & 1.5122033989903259 & z\\ 
720720 &8.560085828505317 & 240 &  1.5618773555849463 & Both\\	
1441440 & 8.907925303844221& 288 & 1.5730156223216971 & Both\\	
2162160	& 9.111778860078964 & 320& 1.5796240997550626 &  Both\\
3603600	& 9.288911580881782	& 360& 1.5781086320000497 & z\\	
4324320	& 9.473895419732957	& 384 & 1.5920793991588504 & Both	\\
7207200	& 9.655498762787248	&  432 & 1.5911044178397928	& z\\
8648640	& 9.680805342657756	& 448 & 1.5857712086389766 & z\\	
10810800 &9.870365907305127	& 480 & 1.598754109907839 & Both \\ 21621600 & 10.251889488064096 & 576	& 1.6129194210253008 & Both\\
36756720 & 10.418517980548113 & 640 & 1.6124072263755087 & z \\
43243200 & 10.469364642977977 &  672 &  1.608133517192604 & z\\	
61261200 & 10.611850914128109 & 720 & 1.6129268471701907 & Both\\	
73513440 &	10.814367000246168 & 768& 1.6277407074294288 & Both\\
122522400 &	11.012526770001736 & 864 & 1.628693078747944 & Both	\\
147026880 & 11.039663838081106 & 896 & 1.6239718455034093 & z\\	
183783600 & 11.247419532881505 & 960 & 1.6379101218780654 & Both\\	
367567200 & 11.664231546026192 & 1152 & 1.6546758215862982 & Both\\	
698377680 & 11.725712741610478 & 1280 & 1.6389019055063054 & z \\	
 \hline
\end{tabular}
\end{center}

    \begin{center}
\begin{tabular}{ |c|c|c|c|c| } 
 \hline
$n$ & $z(n)$ & $\tau(n)$  & $v(n)=\frac{z(n)}{\log \tau(n)}$ & Record type  \\ 
735134400  & 11.90092229151731   &  1344 & 1.6521244370715231 & z\\
1102701600 & 11.928538234423074	& 1440   & 1.6402481809316225 &z\\	
1163962800 & 11.93491248056824  &	1440 & 1.6411246794127394	&z\\
1396755360 & 12.154635463940904	& 1536   & 1.6566362239305026 & Both\\
2205403200 & 12.168527631038101 & 1680   &  1.6385170975662962	& z\\
2327925600 & 12.369007879661227	&1728   & 1.659218315869329 &Both\\
2793510720 & 12.397909801933391	& 1792  & 1.6550213366391964	& z\\
3491888400 & 12.623535748416932	& 1920  &	 1.6697620887531766	& Both \\
6983776800 &	13.074932764428576	& 2304	& 1.6887437164315817 & Both\\
13967553600	&  13.330404681822575	&2688	& 1.688129641564845	& z\\
20951330400	& 13.359790283619853	& 2880	& 1.6771971437444477 & z\\
27935107200	&  13.473027345213657	& 3072 & 1.6778189601814566& z \\
41902660800	&  13.618787163390186	&  3360	& 1.6772532788317212 & z \\
48886437600	& 13.661035931587017 &  3456 & 1.6766395074626765 & z\\
80313433200	&	13.898496037013476	& 3840	& 1.6840073524996981 &z\\
97772875200	& 13.9238684365927 & 4032 & 1.6771667750746457 & z \\
146659312800 &	13.95391871073356 & 4320 &	1.6669335689665656 & z\\
160626866400 &	14.381230476504792 &  4608	& 1.7048362983235954 & Both\\
321253732800 &	14.653665601698984 &	 5376 & 1.70595781024433 & Both\\
481880599200 & 14.684621625013545 & 5760 &1.6959397963111271 & z\\
642507465600 &  14.805417117020358 &  6144 & 1.697240007810422 & z\\	
963761198400 & 14.96078376959392 &  6720 & 1.697611432956445 & z\\
1124388064800 & 15.005991789605918 & 6912 & 1.6973156317382547	& z\\
1927522396800 &  15.114528244214418 & 7680 & 1.6894584161675266 & z\\	
2248776129600 &	15.286212094623194 &	 8064& 1.6993809575575591 & z\\
3373164194400 & 15.317866917412017 & 8640 & 1.6899382326482713 & z \\
4497552259200 & 15.442133591869041 & 9216 & 1.6916033729345097 & z\\	
4658179125600 & 15.532888047165569 & 9216 & 1.701545039464922 & z\\
6746328388800 &  15.601905131290353 & 10080 &  1.6924910964782531 & z\\
9316358251200 &	 15.819921876477741	& 10752 & 1.7042101167364407 & z\\
13974537376800 & 15.85220556834494	& 11520	&  1.695089488250149 & z\\
18632716502400 & 15.979508397400322	& 12288	& 1.6969908838209302 & z\\ 27949074753600 & 16.143138304865357	&  13440 & 1.698206842426234 & z\\
32607253879200 &  16.19090273009808 & 13824 & 1.698198918206643 & z\\	
55898149507200 & 16.304808192777614	&  15360 & 1.6914540141556873 & z\\	
65214507758400 & 16.486083129236103 & 16128	&  1.701646565312295 & z\\
97821761637600 &  16.519094363887366 &	17280 & 1.6929976353281253 & z\\
130429015516800	& 16.65002989232633 & 18432 &  1.695204142143814 & z\\
144403552893600	& 16.681117338585267 & 18432 & 1.698369275660548 & z\\
195643523275200	& 16.818268468291016 &	20160 & 1.6968514958550285 & z\\
288807105787200	& 16.982546580718836 & 21504 & 1.7023412570717529 & z\\
391287046550400	& 16.984348388710064 & 23040 &  1.690828290229381 & z\\
433210658680800	& 17.01612849485301	& 23040	& 1.6939920679206526 & z\\
577614211574400	& 17.14984918244781	& 24576	& 1.6964049353083852 & z\\
866421317361600	& 17.321625217537015 & 26880 &1.69834211149153 & z\\
1010824870255200 & 17.37191487994591 & 27648 & 1.69858126294547	&z\\ 
1732842634723200 & 17.49109979271339 & 30720 &1.6927958764301567 & z\\
2021649740510400 & 17.681845203355937 & 32256 & 1.7032138358195825 &z\\
3032474610765600 & 17.716182715340032 & 34560 & 1.6952551274755545 &z\\
4043299481020800 &	17.85369455966352 &  36864 & 1.6979277485105486 & z\\
\hline
\end{tabular}
\end{center}

    \begin{center}
\begin{tabular}{ |c|c|c|c|c| } 
 \hline
$n$ & $z(n)$ & $\tau(n)$  & $v(n)=\frac{z(n)}{\log \tau(n)}$ & Record type  \\ 
6064949221531200 & 18.030279139683728 & 40320 & 1.700231428280654 & z\\
10685862914126400 & 18.034756898282858 & 43008 &  1.69036627988754 &z\\
12129898443062400 &  18.204352248079783	& 46080 & 1.695299364613963 & z\\
21371725828252800 & 18.208916713887955 & 49152 & 1.6855936510511806	& z\\ 24259796886124800 & 18.299901250390008 & 51840 & 1.6857074987195217	& z\\
30324746107656000 & 18.32631342665929 & 53760 & 1.6825040355068184 & z\\
32057588742379200 &  18.38794334401685 & 53760 & 1.6881621611891613	& z\\
37400520199442400 & 18.44048746990043 & 55296 & 1.6886188378525169& z\\
60649492215312000 & 18.502217106680238	&  61440 & 1.6780813764210956 & z\\
64115177484758400 &  18.5643532977372&  61440 & 1.6837168948248216 & z\\
74801040398884800 & 18.763548743910583 &  64512 & 1.694285814101721	 & z\\
112201560598327200 & 18.799052019527668	 & 69120 & 1.686982036850719& z\\	
149602080797769600 & 18.942419905102103	&  73728 & 1.6900594936024502 & z\\
224403121196654400 & 19.12643161225541 &  80640 & 1.6929416451262043 & z\\	448806242393308800 & 19.30760643717392	&  92160 & 1.6890150485603024 & z\\	
897612484786617600 &  19.406936363910518 & 103680 &  1.6803903332985408	& z\\
\hline
\end{tabular}
\end{center}

\end{document}